\documentclass[10pt]{amsart}

\usepackage[matrix,arrow,curve,frame]{xy}
\usepackage{amsmath,amsthm,amssymb,enumerate}
\usepackage{latexsym}
\usepackage{amscd}

\newtheorem{thm}[subsection]{Theorem}
\newtheorem{defn}[subsection]{Definition}
\newtheorem{prop}[subsection]{Proposition}

\newtheorem{cor}[subsection]{Corollary}
\newtheorem{lemma}[subsection]{Lemma}

\newtheorem{remark}[subsection]{Remark}

\theoremstyle{definition}
\newtheorem{example}[subsection]{Example}

\newcommand{\cat}{\mathcal}
\newcommand{\lra}{\longrightarrow}
\newcommand{\lla}{\longleftarrow}
\newcommand{\llra}[1]{\stackrel{#1}{\lra}}
\newcommand{\llla}[1]{\stackrel{#1}{\lla}}

\newcommand{\R}{\mathbb R}
\newcommand{\Q}{\mathbb Q}
\newcommand{\Z}{\mathbb Z}

\DeclareMathOperator{\id}{id}

\DeclareMathOperator{\colim}{colim}
\DeclareMathOperator{\hocolim}{hocolim}

\DeclareMathOperator{\Met}{Met}

\begin{document}

\title{Homotopy decomposition of a group of 
symplectomorphisms of $S^2\times S^2$}
\author{S\'\i lvia Anjos}\thanks{Supported in part by FCT through program POCTI and grant 
POCTI/1999/MAT/33081.}
\author{Gustavo Granja}\thanks{Supported in part by FCT through program POCTI and grant 
POCTI/1999/MAT/34015.}
\address{Departamento de Matem\'atica \\ 
Instituto Superior T\'ecnico \\  
Av. Rovisco Pais \\
 1049-001 Lisboa\\
 Portugal\\
Fax: (351)218417598}
\date{March 7, 2003}
\email{ggranja@math.ist.utl.pt}
\email{sanjos@math.ist.utl.pt}

{\abstract We continue the analysis started by Abreu, McDuff and Anjos \cite{Ab,AM,McD,An} 
of the topology of the group of symplectomorphisms of $S^2\times S^2$ when 
the ratio of the area of the two spheres lies in the interval $(1,2]$. We express
the group, up to homotopy, as the pushout (or amalgam) of certain of its compact Lie 
subgroups. We use this to compute the homotopy type of the classifying space of the 
group of symplectomorphisms and the corresponding ring of characteristic 
classes for symplectic fibrations.}


\maketitle

\section{Introduction}

Let $M_{\lambda}$ denote the symplectic manifold $(S^2 \times S^2, \omega_\lambda = \lambda\sigma_0 
\oplus \sigma_0 )$ where $1 \leq \lambda \in \R$ and $\sigma_0$ is the standard area form on $S^2$
 with total area equal to 1. It is known that any symplectic form on the manifold $S^2 \times S^2$ is,
 up to scaling by a constant, diffeomorphic to $\omega_\lambda$ (see Lalonde-McDuff \cite{LM}).
 Let $G_\lambda$ denote the group of symplectomorphisms of $M_\lambda$.

By now, much is known about the topology of the group $G_\lambda$ and, more generally, 
about symplectomorphism groups of ruled surfaces. Gromov \cite{Gr} first showed that, when 
$\lambda=1$,  $G_\lambda$ deformation retracts onto the subgroup of standard isometries
$\Z/2 \ltimes ( SO(3)\times SO(3))$. He also showed that this would no longer be the case
if $\lambda >1$, and, confirming this, McDuff \cite{McD1} constructed an element of infinite order in
$H_1(G_\lambda)$. Later, Abreu \cite{Ab} computed the rational cohomology of 
$G_\lambda$ when $1<\lambda\leq 2$ and his methods were extended by Abreu and McDuff \cite{AM,McD}, who
completely described the rational homotopy type of $G_\lambda$ and $BG_\lambda$ for all values of $\lambda$
(see also \cite{McD2} for more information on the integral homotopy 
type of $G_\lambda$, including the fact that the topology of the group changes precisely
 when $\lambda$ crosses an integer). In \cite{An}, the first author computed the 
homology Pontryagin ring of $G_\lambda$ with field coefficients when $1<\lambda \leq 2$ and used
 this to determine the homotopy type of the space $G_\lambda$. All these results were derived from 
a study of the action of $G_\lambda$ on the contractible space $\cat J_\lambda$ of 
compatible almost complex structures on $M_\lambda$. 

Our aim in this paper is to continue the analysis of the case $1<\lambda \leq 2$ by describing 
the homotopy type of $G_\lambda$ as a topological group. 

To explain our result, recall from 
\cite{Ab} that, for $1<\lambda \leq 2$, the space $\cat J_\lambda$ is a stratified space with
 two strata $U_0$ and $U_1$. Each of 
these contains an integrable almost complex structure $J_i$, with isotropy group $K_i \subset G_\lambda$
where $K_0=SO(3)\times SO(3)$ and $K_1=S^1\times SO(3)$. Following Gromov, Abreu showed that 
the inclusions of the orbits $G_\lambda/K_i \subset U_i$ are weak equivalences. $K_0$ and $K_1$ intersect 
in $K_{01} = SO(3)$ which sits inside $K_0$ as the diagonal and inside $K_1$ as the second
factor. Thus, there is a commutative diagram of group monomorphisms
\begin{equation}
\label{square}
\xymatrix{
SO(3) \ar[r]^-{\Delta} \ar[d]_{j} & SO(3)\times SO(3) \ar[d]^{i_0}\\
S^1 \times SO(3) \ar[r]_{i_1} & G_\lambda 
}
\end{equation}
where $\Delta$ denotes the inclusion of the diagonal, $j$ the inclusion of the second 
factor, and $i_0,i_1$ the inclusions in $G_\lambda$ of the  subgroups $K_0,K_1$.

Let $P$ denote the pushout of the diagram of topological groups
\begin{equation}
\label{G-diag}
\xymatrix{
SO(3) \ar[r]^-{\Delta} \ar[d]_{j} & SO(3)\times SO(3) \\
S^1 \times SO(3) & 
}
\end{equation}
This is also known as the amalgam (or amalgamated product, or free product with amalgamation)
 of the groups $SO(3)\times SO(3)$ and $S^1\times SO(3)$ over the common subgroup 
$SO(3)$ and it is characterized as the 
initial topological group admitting compatible homomorphisms from the diagram \eqref{G-diag}.
We will review this construction in section 2. The universal property of the pushout gives us a 
canonical continuous homomorphism 
\begin{equation}
\label{weqpush}
P \lra G_\lambda.
\end{equation}

We need to describe another map derived from \eqref{square}. We will write 
\[ \hocolim (X \llla{f} Y \llra{g} Z) \]
for the homotopy pushout (or double mapping cylinder) of the maps $f,g$. This is the quotient space of 
\[ X \amalg ( Y\times [0,1] )  \amalg Z\]
by the equivalence relation generated by $(y,0) \sim f(y); (y,1) \sim g(y)$.
Applying the classifying space functor to \eqref{square} we 
obtain a canonical map 
\begin{equation}
\label{weqclass}
\hocolim \left( BSO(3)\times BSO(3) \llla{B\Delta}  BSO(3) \llra{Bj} BS^1\times BSO(3) \right) 
\to BG_\lambda.
\end{equation}

Finally, let $\Met(S^2\times S^2)$ denote the space of metrics on $S^2\times S^2$. The 
usual retraction \cite[Proposition 2.50 (ii)]{MS} 
\[ \Met(S^2\times S^2) \llra{r} \cat J_\lambda \]
is equivariant with respect to the action of the group of 
symplectomorphisms. Therefore, letting $K < G_\lambda$ denote a compact subgroup,
 the fact that the fixed  point space $\Met(S^2\times S^2)^K$ is convex
implies that $\cat J_\lambda^K$ is contractible. It follows that we can pick a path 
\begin{equation}
\label{path}
 [0,1] \llra{\gamma} \cat J_\lambda^{K_{01}} 
\end{equation}
with $\gamma(0)=J_0$ and $\gamma(1)=J_1$ and this is unique up to homotopy.
A choice of path determines a $G_\lambda$-equivariant map (unique up to $G_\lambda$-equivariant
homotopy)
\begin{equation}
\label{hocomplex}
\hocolim \left( G_\lambda/(SO(3)\times SO(3)) \llla{\pi_0} G_\lambda/SO(3) \llra{\pi_1} 
G_\lambda/(S^1\times SO(3) \right) \to \cat J_\lambda 
\end{equation}
where $\pi_0,\pi_1$ denote the canonical projections.

We can now state our main result:
\begin{thm}
\label{main} 
Let $1<\lambda \leq 2$. Then the following equivalent statements hold:
\begin{enumerate}[(i)]
\item The homomorphism \eqref{weqpush} is a weak equivalence.
\item The map \eqref{weqclass} is a weak equivalence.
\item The $G_\lambda$-equivariant map \eqref{hocomplex} is a weak equivalence.
\end{enumerate}
\end{thm}

Statement (i) is an immediate consequence of the observation that the homology computations in \cite{An} 
say that, for $k$ a field, the Pontryagin ring $H_\ast(G_\lambda;k)$ is the pushout in the
category of $k$-algebras of the diagram obtained by applying $H_\ast(-;k)$ to  
\eqref{G-diag}, together with Theorem \ref{main1} below which says that the functor $H_\ast(-:k)$
preserves certain pushouts. The equivalence of statements (i)-(iii) is an easy consequence
of Theorem \ref{main1} together with Theorem \ref{main2} which says that under the same 
hypotheses, the classifying functor also preserves pushouts.
From (ii), it is easy to compute the ring $H^*(BG_\lambda;\Z)$ of integral characteristic 
classes of symplectic fibrations with fibre $M_\lambda$ (see Corollary \ref{homology}).

A better proof of Theorem \ref{main}  (which would, in particular, provide an alternative proof of
the main results of \cite{An}) would be to deduce (iii) directly from an analysis of the 
stratification of $\cat J_\lambda$. 
Unfortunately, we were unable to do this. The description of the stratification of 
$J_\lambda$ by McDuff in \cite{McD} yields a homotopy pushout decomposition of $\cat J_\lambda$ as 
\[\hocolim (U_0 \leftarrow NU_1\setminus U_1 \rightarrow NU_1 )\]
 where $NU_1$ denotes a tubular neighborhood of $U_1$ in $\cat J_\lambda$ which fibers over $U_1$
 as a disk bundle. These three spaces have the required homotopy types but we do not know whether 
$NU_1$ can be chosen so as to be invariant under the action of $G_\lambda$. More precisely,
(iii) would follow if we could choose a tubular neighborhood $NU_1$, a path 
$\gamma$ as in \eqref{path}, and $t_0 \in ]0,1[$ satisfying:
\begin{itemize}
\item $G_\lambda \cdot \gamma(t) \subset NU_1$ for $t>t_0$,
\item $\gamma(t) \not\in U_1$ for $t<t_0$.
\end{itemize}

(ii) and (iii) are the statements one should try to generalize 
to the cases when $\lambda>2$, although in those cases it might be necessary 
to use a topological category to index the homotopy colimit decomposition.
See \cite[Appendix D]{MW} for a relationship between stratifications and homotopy colimit 
decompositions.

It is an easy consequence of \cite[Theorem 8, p. 36]{Se} that any compact subgroup of the
amalgam $P$ is subconjugate in $P$ to either $K_0$ or $K_1$. In view of Theorem \ref{main} (ii),
it is natural to ask whether this is also the case in $G_\lambda$. Yael Karshon has proved in \cite{Ka} 
that the answer is yes if one assumes that the subgroup in question is a torus (her result is not
circumscribed to the case $1<\lambda \leq 2$). 

Finally, Theorem \ref{main} suggests it might be a good idea to explore analogies with infinite 
dimensional algebraic groups (cf. \cite{Mi,Ki}).
In particular, compare Theorem \ref{main}(ii) with the homotopy decomposition of the classifying space of a 
rank $2$ Kac-Moody group of indefinite type as the homotopy pushout of a diagram of compact subgroups 
\cite[Theorem 4.2.3]{Ki} (see also \cite{ABKS,BrK}).

\subsection{Organization of the paper} In section 2 we begin by recalling the constructions of pushouts
in various categories and introducing necessary notation. In section 3, under the assumption
that the homomorphisms involved induce injections on homology, we compute the Pontryagin ring 
and the homotopy type of the classifying space of an amalgam of topological groups whose 
connected components are compact Lie groups. In section 4, we deduce Theorem \ref{main} from
 the results of section 3 and \cite{An}. We then use it to produce an interesting fiber sequence
 involving $BG_\lambda$ and compute the ring $H^*(BG_{\lambda};\Z)$ of characteristic classes for
 symplectic fibrations with fiber $M_\lambda$.

\section{Pushouts and amalgams.}

In this section we begin by reviewing the categorical notions of pushout and coequalizer.
We then recall the construction of the pushout in the categories of 
groups, topological groups and $k$-algebras, and introduce some notation that will be 
necessary in the next section.

\subsection{Pushouts and coequalizers}
Everything in this subsection is standard basic category theory. See \cite{Ma}, for instance.

Let $B_0, B_1, B_2$ be objects in a category $\cat C$
 and $f_1,f_2$ morphisms in $\cat C$. Given a diagram in $\cat C$,
\begin{equation}
\label{pushout}
\xymatrix{
B_0 \ar[d]^{f_2} \ar[r]^{f_1} & B_1 \\
B_2 & \\
}
\end{equation}
a \emph{cone} on this diagram consists of an object $C$ together with morphisms $g_1$ and $g_2$ such
 that
$$\xymatrix{
B_0 \ar[d]^{f_2} \ar[r]^{f_1} & B_1  \ar[d]^{g_1}\\
B_2 \ar[r]^{g_2} & C \\
}
$$
commutes. A \emph{pushout} of (\ref{pushout}) is a cone $(C,g_1,g_2)$ with the property that 
if $(D,h_1,h_2)$ is any other cone, then there is a unique morphism
 $\varphi : C \rightarrow D$  making the following diagram commute: 
$$\xymatrix{
B_0 \ar[d]^{f_2} \ar[r]^{f_1} & B_1  \ar[d]^{g_1} \ar[ddr]^{h_1} &\\
B_2 \ar[r]^{g_2} \ar[drr]^{h_2} & C   \ar@{-->}[dr]^{\varphi}  & \\
  &  & D }
$$
This universal property characterizes the pushout up to unique isomorphism in $\cat C$
(when the pushout exists). The pushout is usually denoted by
\[ B_1\amalg_{B_0} B_2. \]

Given a pair of arrows
\begin{equation}
\label{coequalizer}
\xymatrix{ A \ar@<1ex>[r]^{d_0} \ar@<-1ex>[r]_{d_1} & B},
\end{equation}
a cone on \eqref{coequalizer} is an arrow $\epsilon:B \to C$ such that $\epsilon d_0= \epsilon d_1$. 
A \emph{coequalizer} of \eqref{coequalizer} is a cone $\epsilon:B \to C$ with the property 
that given any other cone  $\epsilon':B \to C'$, there is a unique map $\varphi: C \to C'$ making
the following diagram commute 
\[
\xymatrix{A \ar@<1ex>[r]^{d_0} \ar@<-1ex>[r]_{d_1} & B \ar[r]^\epsilon \ar[dr]^{\epsilon'}& C 
\ar[d]^{\varphi} \\
& & C'}
\]
Again this universal property characterizes the coequalizer up to unique isomorphism.
Pushouts and coequalizers are instances of a more general construction, that of 
colimit of a diagram, which is defined by a similar universal property.

\begin{example} 
\label{example}
The main examples we have in mind are the following:
\begin{enumerate}[(a)]
\item $\cat C$ is the category of topological spaces. 
If $A$ is a topological group acting on the right on the space $Y$ and on the left on the 
space $X$ then the coequalizer of the two maps $d_0,d_1:X\times A \times Y \to X\times Y$ 
defined by $d_0(x,a,y) = (xa,y)$ and $d_1(x,a,y)=(x,ay)$ is the quotient of $X\times Y$
by the action $(a,x,y) \to (xa,a^{-1}y)$. We write $X\times_A Y$ for this space.  
\item $\cat C$ is the category of vector spaces over a field $k$. If $A$ a $k$-algebra, $V$
 a left $A$-module and $W$ a right $A$-module,
the coequalizer of the action maps $d_0,d_1 : V \otimes_k A \otimes_k W \to V \otimes_k W$
defined as above is the tensor product $V \otimes_A W$.
\end{enumerate}
\end{example}

If $k$ is a field, the universal property of the coequalizer together with the
K\"unneth theorem give us a canonical map 
\begin{equation}
\label{examplekun}
 H_*(X;k) \otimes_{H_*(A;k)} H_*(Y;k) \to H_*(X \times_A Y ; k). 
\end{equation}
The main theorem of this paper follows immediately from the fact that a similar canonical map is an 
isomorphism and the proof will consistently exploit such maps.

\subsection{Pushouts of groups}

Suppose $\cat C$ is the category of groups. A good reference for everything in 
this subsection is \cite{Se}.

Let $S$ denote the set of finite sequences 
\[ x_1 x_2 \ldots x_n \]
where $x_i \in B_1$ or $x_i \in B_2$ and consider the equivalence relation
$\sim$ on $S$ generated by
\begin{enumerate}[(i)]
\item $x_1\ldots x_n \sim x_1\ldots \hat{x_i} \ldots x_n $ if $x_i=1$,
\item $x_1\ldots f_1(a) \ldots x_n  \sim x_1 \ldots f_2(a) \ldots x_n$ for $a \in B_0$,
\item $x_1 \ldots x_i x_{i+1} \ldots x_n \sim x_1 \ldots (x_i x_{i+1}) \ldots x_n$ when $x_i,x_{i+1}$
both belong to $B_1$ or $B_2$.
\end{enumerate}
$S$ has an associative unital product defined by concatenation and it is easy to check that this
descends to the set $P=S/\sim$ of equivalence classes and that $P$ together with the canonical
maps $B_1 \to P$ and $B_2 \to P$ is the pushout of \eqref{pushout}.

If the maps $f_1$ and $f_2$ are monomorphisms 
then the pushout of the diagram \eqref{pushout} is also called the 
\emph{amalgam}\footnote{It is also called amalgamated product of $B_1$ and $B_2$ over $B_0$, 
and free product of $B_1$ and $B_2$ with amalgamation.} of $B_1$ and $B_2$ over $B_0$.
In this case, there is a useful description of the elements of the pushout, called 
the Normal Form Theorem. To explain this, we will start by introducing some notation.
Let ${\bf i}=(i_1,\ldots,i_n) \in \{0,1,2\}^{n}$. Then we define
\begin{equation}
\label{defBis}
\overline{B}_{\bf i} = B_{i_1} \times \cdots \times B_{i_n} 
\end{equation}
Note that there is a canonical map 
\begin{equation}
\label{projection}
\overline{B}_{\bf i} \llra{\pi} P 
\end{equation}
determined by multiplication. $B_0^{n-1}$ acts on the right on $\overline {B}_{\bf i}$ by 
\[ (a_1,\ldots,a_{n-1}) \cdot (b_1,\ldots,b_n) = (b_1a_1, a_1^{-1}b_1a_2,a_2^{-1}b_2a_3,\ldots,
a_{n-1}^{-1} b_{n}). \]
We will write 
\begin{equation}
\label{defBis2}
B_{\bf i} = B_{i_1} \times_{B_0} \cdots \times_{B_0} B_{i_n} 
\end{equation}
for the quotient of $\overline{B}_{\bf i}$ by this action. Note that we can express the 
quotient $B_{\bf i}$ as the coequalizer of 
\[ \xymatrix{ \overline{B}_{\bf i} \times B_0^{n-1} \ar@<1ex>[r]^-{d_0} \ar[r]_-{d_1} &
\overline{B}_{\bf i}
}\]
where $d_0$ is given by the action, and $d_1$ is the projection on the first factor.
The canonical map $\pi$ of \eqref{projection} factors through this quotient and we will 
denote the resulting map also by $\pi$.

Let $B_i'=B_i \setminus B_0$ for $i=1,2$, and for $\mathbf{i}=(i_1,\ldots,i_n)$ a sequence of
alternating $1$'s and $2$'s, let
\[ B_{\mathbf{i}}' = B_{i_1}' \times_{B_0} B_{i_2}' \times_{B_0} \cdots \times_{B_0} B_{i_n}' \]
denote the quotient of $B_{i_1}' \times \cdots \times B_{i_n}'$ by the action of 
$B_0^{n-1}$.

\begin{thm}\cite[Theorem 2, page 4]{Se}
\label{normalform}
If $f_1, f_2$ are monomorphisms, the canonical maps $\pi$ determine a bijection 
\[ B_0 \amalg \amalg_{\bf i} B_{\bf i}' \llra{\sim} P \]
where ${\bf i}$ runs over all sequences of alternating $1$'s and $2$'s.
\end{thm}

\subsection{Pushouts of topological groups}

We now consider the case when the $B_i$ are topological groups and the homomorphisms between them 
continuous. We need to give a topology to the space $P$ defined above so that 
$P$ together with the maps $B_i \to P$ is the pushout in the category of topological groups. 

In order to do this, it is necessary to work in a 
category of compactly generated spaces. Any will do for our purposes, so we will work 
with the simplest, namely that of Vogt \cite[Example 5.1]{Vo}. A space $X$ is said to be 
compactly generated if a subset $U \subset X$ is closed iff for every compact Hausdorff space 
$K$ and continuous map $g:K \to X$, $g^{-1}(U)$ is closed. Given an arbitrary
topological space $X$ we can refine the topology in the obvious way so that it
becomes compactly generated. Denoting this space by $kX$, the natural transformation
\[ kX \to X \]
determined by the identity maps is a weak equivalence \cite[Proposition 1.2 (h)]{Vo}. 
The (categorical) product in the category of compactly generated spaces does
not in general agree with the product in spaces; it is necessary to apply the functor $k$ to the usual
product. In all that follows we will write $\times$ for the product in the category of 
compactly generated spaces. By a topological group we mean a group object in the category of
compactly generated spaces. Any topological group in the usual sense
determines such a group object by applying $k$ to the multiplication\footnote{This is not 
necessary for groups of diffeomorphisms of compact manifolds which are metrizable and
hence compactly generated.}.

We will write 
\[ P_n \subset P \]
for the image in $P$ of the $B_{\bf i}$ with ${\bf i}\in \{0,1,2\}^n$ and $P_0$ for the image 
of $B_0$. We will also write ${\bf \alpha}_n$ and ${\bf \beta}_n$ for the two sequences
in $\{0,1,2\}^n$ of alternating $1$'s and $2$'s and 
\begin{equation}
\label{defQn}
Q_n = B_{{\bf \alpha}_n} \amalg B_{{\bf \beta}_n}. 
\end{equation}
Similarly we define 
\begin{equation}
\label{defQn2}
\overline{Q}_n = \overline{B}_{{\bf \alpha}_n} \amalg \overline{B}_{{\bf \beta}_n}. 
\end{equation}
We will set $\overline{Q}_0 = Q_0 = B_0$.
With the notation above, we first give $P_n$ the quotient topology determined by the canonical map
\begin{equation}
\label{projn}
Q_n \llra{\pi} P_n.
\end{equation}
Note that this is the same as to give $P_n$ the topology induced by $\overline{Q}_n$, since
$Q_n$ has the quotient topology determined by the projection $\overline{Q}_n \to Q_n$.

The inclusion  $\overline{Q}_n \to \overline{Q}_{n+1}$ defined by 
\[ (x_1,\ldots,x_n) \mapsto (x_1,\ldots,x_n,1) \]
induces a continuous map 
\[ P_n \llra{i_n} P_{n+1}. \]
We give $P$ the topology of the union of the $P_n$'s. The product $\mu:P\times P \to P$ clearly restricts 
to a product $\mu_{n,m}: P_n \times P_m \to P_{n+m}$.
Consider the commutative diagram
\[
\xymatrix{
\overline{Q}_n \times \overline{Q}_m \ar[r]^{\overline{\mu}_{n,m}} \ar[d]_{\pi \times \pi} & 
\overline{Q}_{n+m} \ar[d]^\pi \\
P_n \times P_m \ar[r]^{\mu_{n,m}}&  P_{n+m} 
}\]
where the map $\overline{\mu}_{n,m}$ is defined by:
\[
\overline{\mu}_{n,m} (\mathbf{x},\mathbf{y})  =   \left\{ \begin{array}{ll}
(x_1,\ldots,x_n,y_1,\ldots,y_m) & \text{ if } x_n \in B_i, y_1 \in B_j \text{ with } i \neq j,\\
(x_1,\ldots,x_n y_1, \ldots, y_m, 1) & \text{ otherwise.}
\end{array}
\right.
\]
Both vertical maps are quotient maps \cite[Corollary 3.8]{Vo} so we conclude that 
the maps $\mu_{n,m}$ are continuous.
Since in the category of compactly generated spaces, colimits commute with products
\cite[Proposition 3.7 (b)]{Vo}, it follows that the multiplication
\[ P\times P \llra{\mu} P \]
is continuous.
Similarly one checks that the inverse map on $P$ is continuous. It is now easy to check that the 
topological group $P$ together with the maps $B_i \to P$ has the required 
universal property and hence is the pushout in the category of topological groups.

\subsection{Pushouts of $k$-algebras}

Let $k$ be a field. There is an entirely analogous description of the pushout of \eqref{pushout}
 in the category of associative $k$-algebras. We will soon want to compare the construction
for algebras and the homology of the analogous construction for topological groups and 
to make the notation more clear we will decorate the algebraic constructions with a superscript \emph{alg}. 

Let $S^{alg}$ denote the tensor algebra generated 
by the $k$-vector spaces $B_1^{alg}$ and $B_2^{alg}$. Thus, as a vector space, $S^{alg}$ is the direct sum
\[ 
k \oplus B_1^{alg} \oplus B_2^{alg} \oplus (B_1^{alg} \otimes B_1^{alg}) \oplus (B_1^{alg}\otimes B_2^{alg})
 \oplus \ldots 
\]
where $\otimes$ denotes the tensor product of $k$-vector spaces.
Then the pushout of \eqref{pushout} in $k$-algebras is 
\[ P^{alg}=S^{alg}/I \]
together with the canonical maps $B_i^{alg} \to P^{alg}$, where $I$ is the ideal generated by 
\begin{enumerate}[(i)]
\item $x_1\otimes\ldots \otimes x_n - x_1\otimes \ldots \otimes \hat{x_i} \otimes \ldots \otimes 
x_n $ if $x_i=1$,
\item $x_1\otimes \ldots \otimes f_1(a) \otimes \ldots x_n - x_1 \otimes f_2(a) \otimes \ldots x_n$ 
for $a \in B_0^{alg}$,
\item $x_1 \otimes \ldots \otimes x_i \otimes x_{i+1} \otimes \ldots \otimes x_n - x_1 \otimes
 \ldots\otimes  (x_i x_{i+1}) \otimes \ldots \otimes x_n$ when $x_i,x_{i+1}$
both belong to $B_1^{alg}$ or $B_2^{alg}$.
\end{enumerate}

In some circumstances it is possible to give a more precise description 
of the pushout. To do this, we need some notation.
Let ${\bf i}=(i_1,\ldots,i_n) \in \{0,1,2\}^{n}$. Then we define 
\begin{equation}
\label{defBisalg}
\overline{B}^{alg}_{\bf i} = B_{i_1}^{alg} \otimes \cdots \otimes B_{i_n}^{alg} 
\end{equation}
Note that there is a canonical map 
\begin{equation}
\label{projectionalg}
\overline{B}^{alg}_{\bf i} \llra{\pi} P^{alg} 
\end{equation}
determined by the multiplication on $P^{alg}$. 

Given a sequence ${\bf i} \in \{0,1,2\}^n$, let $\overline{\bf i} \in \{0,1,2\}^{2n-1}$ denote
the sequence obtained from ${\bf i}$ by inserting $0$'s between each two entries of ${\bf i}$.
For example, if ${\bf i}= (i_1,i_2,i_3)$ then $\overline{\bf i} = (i_1,0,i_2,0,i_3)$.
Then we define maps
\[ \xymatrix{ \overline{B}^{alg}_{\overline{\bf i}} \ar@<1ex>[r]^{d_0} \ar[r]_{d_1} &
\overline{B}^{alg}_{\bf i}
}\]
by (setting $f_0 = \id$ )
\begin{align*}
d_0(x_1\otimes a_1\otimes x_2 \otimes \ldots \otimes x_n) & = x_1 f_{i_1}(a_1) \otimes x_2 f_{i_2}(a_2)
 \otimes \ldots \otimes x_n \\
d_1(x_1\otimes a_1 \otimes x_2 \otimes \ldots \otimes  x_n) & = x_1 \otimes f_{i_2}(a_1) x_2 \otimes 
\ldots \otimes f_{i_n}(a_{n-1}) x_n.
\end{align*}
and we define $B_{\bf i}^{alg}$ to be the coequalizer (in $k$-vector spaces) of these maps. That is,
\begin{equation}
\label{defBis2alg}
B_{\bf i}^{alg} = B_{i_1}^{alg} \otimes_{B_0^{alg}} B_{i_1}^{alg} \otimes_{B_0^{alg}} \ldots 
\otimes_{B_0^{alg}} B_{i_n}^{alg}
\end{equation}
The canonical map $\pi$ of \eqref{projectionalg} factors through $B_{\bf i}^{alg}$ and we will 
denote the resulting map also by $\pi$. Moreover we will write as in \eqref{defQn}
\begin{equation}
\label{defQnalg}
Q_n^{alg} = B_{{\bf \alpha}_n}^{alg} \oplus B_{{\bf \beta}_n}^{alg}. 
\end{equation} 
and
\begin{equation}
\label{defPnalg}
P_n^{alg} \subset P^{alg}
\end{equation}
for the image of $Q_n^{alg}$ in $P^{alg}$.

\begin{remark}
\label{equivdefalg}
If the $k$-algebras $B_i^{alg}$ are the homology Pontryagin rings of topological groups $B_i$, 
then they have added structure: they are Hopf algebras (see \cite{MM}). In particular there 
are antiautomorphisms 
\[ c:B_i^{alg} \to B_i^{alg}\]
induced by the inverse map on the group, as well as augmentations (i.e. maps of $k$-algebras)
\[ B_i^{alg} \llra{\epsilon} k \]
induced by the map to the trivial group. Let $B_0^{n \ alg}$ denote the $k$-algebra
\[ B_0^{n \ alg}  = B_0^{alg} \otimes \cdots \otimes B_0^{alg} \]
augmented in the obvious way. For  ${\bf i}=(i_1,\ldots,i_n) \in \{0,1,2\}^n$, 
 $\overline{B}_{\bf i}^{alg}$ is a right $B_0^{n-1 \ alg}$-module under the action 
\[ (b_1\otimes \ldots \otimes b_n) \cdot (a_1\otimes\ldots\otimes a_{n-1})= (b_1a_1 \otimes c(a_1)b_1a_2
\otimes \ldots \otimes c(a_{n-1}) b_{n}) \]
where we have omitted the maps $f_i$ from the notation. Moreover, it is easy to see that 
\[ \overline{B}_{\bf i}^{alg} \otimes_{B_0^{n-1\ alg}} k = B_{\bf i}^{alg}. \]
As in \eqref{examplekun} it is then an immediate consequence of the K\"unneth theorem and
the universal property of a coequalizer that there is a canonical map
\begin{equation}
\label{canonicalmap}
 B_{\bf i}^{alg} \llra{\phi} H_*(B_{\bf i};k)
\end{equation}
Moreover, writing $P^{alg}$ for the pushout of the Pontryagin rings, the diagram
\[ \xymatrix{ 
B_{\bf i}^{alg} \ar[r]^-{\phi} \ar[d] & H_*(B_{\bf i};k) \ar[d] \\
P^{alg} \ar[r] & H_*(P;k) 
}\]
commutes. 
\end{remark}

We will need the following analog of Theorem \ref{normalform} for pushouts of $k$-algebras, due to 
P.M. Cohn.
\begin{thm}\cite[Proof of Theorem 3.1]{Co3}
\label{normalformalg}
If $f_1, f_2$ are monomorphisms, and there exist right $B_0^{alg}$-modules $B_i'^{alg}$ such that 
$B_i^{alg}= f_i(B_0^{alg}) \oplus B_i'^{alg}$ as right $B_0^{alg}$-modules. Then
regarding $B_i'^{alg}$ as a $B_0^{alg}$-bimodule via the isomorphism $B_i'^{alg} \simeq B_i^{alg}/B_0^{alg}$,
the canonical maps $\pi$ determine an isomorphism of $B_0^{alg}$-bimodules
\[ B_0^{alg} \oplus \oplus_{\bf i} B_{\bf i}'^{alg}  \llra{\sim} \oplus_{n\geq 0} P_n^{alg}/P_{n-1}^{alg} \]
where ${\bf i}$ runs over all sequences of alternating $1$'s and $2$'s.
\end{thm}

For future reference, we also note that there are exhaustive filtrations 
\[ B_i^{alg} \subset \ldots \subset P_n^{alg}B_i^{alg} \subset \ldots \subset P^{alg}\]
of $P^{alg}$ by right $B_i^{alg}$-modules and that (cf. \cite[(32), p. 61]{Co3})  we have 
canonical isomorphisms
\begin{equation}
\label{secondfiltr}
 \oplus_{\bf i} B_{\bf i}'^{alg} \otimes_{B_0^{alg}} B_j^{alg} \llra{\sim} \oplus_{n \geq 0}
P_n^{alg}B_j^{alg}/P_{n-1}^{alg}B_j^{alg} 
\end{equation}
of right $B_j^{alg}$-modules, where ${\bf i}$ runs over all sequences of alternating $1$'s 
and $2$'s not ending in $j$.

\begin{remark}
Theorem \ref{normalformalg} also holds under the assumption that the $f_i$ are monomorphisms 
and  $B_i'^{alg}=B_i^{alg}/f_i(B_0^{alg})$ are flat right $B_0^{alg}$-modules 
(cf. \cite[(1.6)-(1.8) p. 436]{Co2})  but the above formulation 
is sufficient for our purposes. 
\end{remark}

\section{Amalgams of topological groups.}

In this section, under the assumption that the homomorphisms involved induce injections on homology,
 we compute the Pontryagin ring of an amalgam of topological groups whose connected components are
 compact Lie groups. Using this and Puppe's theorem, we obtain a homotopy decomposition of the
classifying space.

Although the results are stated for the case when the free products involved have only two factors,
it is clear that all the statements hold (with the same proof) for an 
arbitrary number of factors.

\subsection{Homology of the pushout}
In order to compute the homology of the pushout of a diagram of topological groups
we require some assumptions. 
\begin{defn}
\label{splitness}
We will say that a pushout diagram \eqref{pushout} of topological groups is \emph{homologically free}
if it satisfies the following conditions:
\begin{enumerate}[(i)]
\item The connected components of the $B_i$ are compact Lie groups, 
\item The continuous homomorphisms $f_1, f_2$ are monomorphisms.
\item The homomorphims $H_*(f_i;k)$ are injective for every field $k$.
\end{enumerate}
\end{defn}
The first condition of the previous definition can certainly be weakened. We make this 
assumption because it suffices for the application we have in mind and it considerably
simplifies the point set topology involved.

\begin{lemma}
\label{relhom}
Under the assumptions (i)-(ii) of Definition \ref{splitness}, and with the notation 
of \eqref{projn} we have
\begin{enumerate}[(a)]
\item $P_n$ is a Hausdorff space.
\item $(Q_n,\pi^{-1}(P_{n-1})) \llra{\pi_n} (P_n,P_{n-1})$ is a relative homeomorphism.
\end{enumerate}
\end{lemma}
\begin{proof}
(a) is easy and (b) follows easily from compactness and the Normal Form Theorem \ref{normalform}.
\end{proof}

The following lemma will be used often in the rest of this section.
\begin{lemma}
\label{tensor}
Let $X \to B$ be a principal fibration with fiber $A$ and $k$ be a field. If
the inclusion of the fiber induces an injection $H_*(A;k) \to H_*(X;k)$ then 
\begin{enumerate}[(a)]
\item  $H_*(X;k)$ is a free (and hence flat) right $H_*(A;k)$-module,
\item  The canonical map 
\[ H_*(X;k) \otimes_{H_*(A;k)} k \to H_*(B;k) \]
is an isomorphism.
\end{enumerate}
\end{lemma}
\begin{proof}
Since the inclusion of the fiber $A \to X$ induces an injection on connected components, it
suffices to consider the case when $A$ is connected. It follows that the action of $\pi_1(B)$ on 
$H_*(A;k)$ is trivial. The homology Leray-Serre spectral sequence is a spectral sequence of 
right $H_*(A;k)$-modules. Since $H_*(A;k) \to H_*(X;k)$ is an injection, spectral sequence collapses.  
The $E_2$ and hence the $E_\infty$
terms are free $H_*(A;k)$-modules on $H_*(B;k)$. It follows that $H_*(X;k)$ is a 
free $H_*(A;k)$-module and the projection induces an isomorphism
\[ H_*(X;k) \otimes_{H_*(A;k)} k \to H_*(B;k)  \]
as required.
\end{proof}
\begin{remark} 
The first part of Lemma \ref{tensor} is also a consequence of \cite[Theorem 4.4]{MM}.
\end{remark}

In order to compute the homology of the pushout, we need one more lemma.
Let $\Lambda$ denote the partially ordered set $\{0,1,2\}$ where the order is defined
by $0\leq 1$ and $0\leq 2$ ($1,2$ are incomparable). Given $n>0$, $\Lambda^n$ is again a poset.
 The order relation is 
\[ (i_1, \ldots, i_n) \leq (j_1,\ldots, j_n) \iff i_k \leq j_k 
\text{ for each } i. \]
For ${\bf j} \leq  {\bf i} \in \Lambda^n$  there are obvious inclusions 
\[
B_{\bf j} \lra B_{\bf i}.
\]
\begin{lemma}
\label{poset}
Suppose the diagram \eqref{pushout} satisfies conditions (i) and (ii) of Definition \ref{splitness}.
Let ${\bf i} \in \Lambda^n$ and $\Pi \subset \Lambda^n$ be such that \footnote{ These 
conditions ensure that the colimit is just the union of the images of the spaces $B_{\bf j}$ 
in $B_{\bf i}$.}
\begin{enumerate}[(a)]
\item ${\bf j} \leq {\bf i}$ for every ${\bf j} \in \Pi$, 
\item ${\bf j} \in \Pi$ and ${\bf k} \leq {\bf j} \Rightarrow {\bf k} \in \Pi$.
\end{enumerate}
Then the canonical map
\[ \colim_{{\bf j} \in \Pi} B_{\bf j} \to B_{\bf i} \]
is a closed cofibration.
\end{lemma}
\begin{proof}
If $\Pi=\{{\bf j}\}$ consists of a single element, this follows from the fact that 
$\overline{B}_{\bf j} \subset \overline{B}_{\bf i}$ has a $B_{0}^{n-1}$-equivariant
tubular neighborhood, from which we can obtain a neighborhood deformation retraction
of $B_{\bf j} \subset B_{\bf i}$.

The result now follows by induction using the union theorem for cofibrations
\cite[Corollary 2]{Li} (which states that if $A \subset X$, $B \subset X$ and $A\cap B \subset X$ 
are closed cofibrations then $A \bigcup B \subset X$ is again a closed cofibration).
\end{proof}

For the rest of this section we set
\[ B_i^{alg} := H_*(B_i;k). \]
\begin{cor}
\label{splitcoeq}
If the pushout diagram \eqref{pushout} is homologically free, then for any ${\bf i} \in \Lambda^n$ 
the canonical map
\[ B_{\bf i}^{alg} \llra{\phi} H_*(B_{\bf i};k) \]
is an isomorphism.
\end{cor}
\begin{proof}
It is easy to check that, under our assumptions, the inclusion of the fiber 
in the total space of the principal $B_0^{n-1}$-fibration
\[ \overline{B}_{\bf i} \to B_{\bf i} \]
is an injection. The result now follows from immediately from 
 Lemma \ref{tensor} and Remark \ref{equivdefalg}.
\end{proof}

Note that if the diagram \eqref{pushout} is homologically free, then the inclusions 
$ B_{\bf j} \subset B_{\bf i}$
induce inclusions on homology with field coefficients, i.e. the canonical maps 
\[
B_{\bf j}^{alg} \lra B_{\bf i}^{alg} 
\]
are injective. This follows from the obvious fact that 
\[ \overline{B}_{\bf j}^{alg} \lra \overline{B}_{\bf i}^{alg} \]
is injective, together with the fact that the $H_*(B_0^{n-1};k)$-modules $\overline{B}_{\bf i}^{alg}$
are flat by Lemma \ref{tensor}(a). More generally, the same argument implies that the canonical map 
\begin{equation}
\label{inclusion}
\colim_{{\bf j} \in \Pi} B_{\bf j}^{alg} \lra B_{\bf i}^{alg}
\end{equation}
is injective, when $\Pi \subset \Lambda^n$ satisfies the two conditions of Lemma \ref{poset}.
\begin{lemma}
\label{homolcolim}
Suppose the pushout diagram \eqref{pushout} is homologically free. Let $\Pi \subset \Lambda^n$
satisfy the two conditions of Lemma \ref{poset}. Then the canonical map 
\[ \colim_{{\bf j} \in \Pi} B_{\bf j}^{alg} \lra H_*(\colim_{{\bf j} \in \Pi} B_{\bf j};k) \]
is an isomorphism.
\end{lemma}
\begin{proof}
Suppose $\Pi' = \Pi \bigcup \{{\bf i}\}$ and let $\Gamma = \{ {\bf j} \in \Pi: {\bf j} \leq {\bf i}\}$.
Then we have a pushout square
\begin{equation}
\label{colimdiag}
\xymatrix{
\colim_{{\bf j} \in \Gamma} B_{\bf j} \ar[r] \ar[d]_{\eta} & \colim_{{\bf j} \in \Pi} B_{\bf j} \ar[d] \\
B_{\bf i} \ar[r] & \colim_{{\bf j} \in \Pi'} B_{\bf j}
}
\end{equation}
By Lemma \ref{poset} the map $\eta$ is a cofibration. If we assume that the result is 
true for $\Gamma$ then $\eta$ induces an injection in homology by \eqref{inclusion}
and so the Mayer-Vietoris sequence for \eqref{colimdiag} splits. This implies that 
the result holds for $\Pi'$, and hence it is true in general, by induction.
\end{proof}

\begin{thm}
\label{main1}
If the pushout diagram \eqref{pushout} is homologically free then the canonical map
\[ P^{alg} = H_*(B_1;k) \amalg_{H_*(B_0;k)} H_*(B_2;k) \lra H_*(P;k) \]
is an isomorphism.
\end{thm}
\begin{proof}
Since homology commutes with sequential colimits of $T_1$ spaces along closed inclusions, it suffices to 
show that the canonical map 
\[ P_n^{alg} \llra{\gamma_n} H_*(P_n;k) \]
is an isomorphism for each $n$.
This is obvious for $n=0$. Assume it is true for $n\leq m$ and consider the following
diagram
\[\xymatrix{
0 \ar[r] & P_m^{alg} \ar[r] \ar[d]_{\gamma_m} & P_{m+1}^{alg} \ar[d]_{\gamma_{m+1}} \ar[r] & 
P_{m+1}^{alg}/P_m^{alg} \ar[d]^{\phi_{m+1}} \ar[r] & 0 \\
& H_*(P_m;k) \ar[r] & H_*(P_{m+1};k) \ar[r] & H_*(P_{m+1},P_m;k) &
}\]
It is enough to see that the induced map $\phi_{m+1}$ is an isomorphism since it then
follows that $\psi_n$ is surjective and then, by the 5-lemma, that $\gamma_{m+1}$ is an
isomorphism.
To simplify notation, write 
\[ W_{n} = \pi^{-1}(P_{n-1}) \subset Q_n \quad ; \quad 
W_{n}^{alg} = \pi^{-1}(P_{n-1}^{alg}) \subset Q_n^{alg}. \]
It is clear from the Normal Form Theorem \ref{normalform} that
\[ W_{n} = \colim_{{\bf j}<\alpha_n} B_{\bf j} \amalg \colim_{{\bf j}<\beta_n} B_{\bf j}. \]
Since $B_1^{alg}$ and $B_2^{alg}$ are free right $B_0^{alg}$-modules by 
Lemma \ref{tensor}, it follows from Theorem \ref{normalformalg} that we have 
similarly
\begin{equation}
\label{tricky}
W_{n}^{alg} = \colim_{{\bf j}<\alpha_n} B_{\bf j}^{alg} \oplus \colim_{{\bf j}<\beta_n} 
B_{\bf j}^{alg}.
\end{equation}

Consider the diagram
\[ \xymatrix{
W_{m+1}^{alg} \ar[r]^{\iota} \ar[d] & Q_{m+1}^{alg} \ar[r]  \ar[d] & Q_{m+1}^{alg}/ W_{m+1}^{alg} 
\ar[d]^{\psi_{m+1}} \\
H_*(W_{m+1},k) \ar[r]^{i_{m+1}} & H_*(Q_{m+1}; k) \ar[r] & H_*(Q_{m+1},W_{m+1};k) 
}\]
Lemma \ref{homolcolim} says that the left vertical map is an isomorphism
and Lemma \ref{splitcoeq} says that the middle vertical one is. By \eqref{inclusion}, $\iota$ is
an inclusion so $\psi_{m+1}$ is an isomorphism.

The desired result now follows from the commutativity of the diagram
\[\xymatrix{
Q_{m+1}^{alg}/ W_{m+1}^{alg}  \ar[r]^{\eta}
 \ar[d]^{\psi_{m+1}} & P_{m+1}^{alg}/P_{m}^{alg} \ar[d]^{\phi_{m+1}} \\
H_{*}(Q_{m+1},W_{m+1};k) \ar[r]^{\xi} & H_*(P_{m+1},P_m;k) }\]
where $\eta$ is an isomorphism by definition and $\xi$ by Lemma \ref{relhom}(b).
\end{proof}

\subsection{The homotopy decomposition of $BP$.}

We will need the following corollary of Puppe's theorem (see \cite[Proposition, page 180]{DF} or
\cite[Theorem 1]{MaWa}).
\begin{thm}
\label{puppe}
Given a diagram of spaces 
\[ \xymatrix{ 
X_1 \ar[d] & \ar[l] X_0 \ar[r] \ar[d] & X_2 \ar[d] \\
Y & \ar[l]_{=} Y \ar[r]^{=} & Y 
}
\]
then denoting by $F_i$ the homotopy fiber of the map $X_i \to Y$, and
by $F$ the homotopy fiber of the canonical map 
\[ \hocolim (X_1 \lla X_0 \lra X_2) \to Y, \]
the canonical map 
\[ \hocolim (F_1 \lla F_0 \lra F_2) \to F \]
is a weak equivalence.
\end{thm}
We now compute the homotopy type of the classifying space $BP$. Note that this result generalizes the
 well known theorem of J.H.C. Whitehead for discrete groups \cite[Theorem II.7.3]{Br}. 
\begin{thm}
\label{main2}
If the pushout diagram \eqref{pushout} is homologically free then, the canonical map 
\[ \hocolim( BB_1 \lla BB_0 \lra BB_2) \llra{\phi} BP \]
is a weak equivalence.
\end{thm}
\begin{proof}
Consider the diagram 
\[ \xymatrix{
P/B_1 \ar[d] & \ar[l] P/B_0 \ar[d]\ar[r] & P/B_2 \ar[d]\\
BB_1 \ar[d] & \ar[l] BB_0 \ar[d]\ar[r] & BB_2 \ar[d] \\
BP & \ar[l]_{=} BP \ar[r]_{=} & BP
}\]
where all the vertical maps are fiber sequences.
It follows from Van-Kampen's theorem that $\phi$ induces an isomorphism on $\pi_1$.
Therefore the homotopy fiber of $\phi$, which by Theorem \ref{puppe} is
\begin{equation}
\label{eq1}
F \simeq \hocolim (P/B_1 \lla P/B_0 \lra P/B_2),
\end{equation}
 has abelian  fundamental group. It now suffices to show that for any field $k$, we have 
\begin{equation}
\label{eq2}
H_\ast(F;k)\simeq k.
\end{equation}
By Lemma \ref{tensor} it is sufficient to show that the map
\begin{equation}
\label{bruhat}
P^{alg}\otimes_{B_0^{alg}} k \to P^{alg}\otimes_{B_1^{alg}} k \oplus 
P^{alg}\otimes_{B_2^{alg}} k
\end{equation}
is an isomorphism in positive degrees and injective with cokernel $k$ in degree 0.

Consider the filtrations $P_n^{alg} , P_n^{alg}B_1^{alg}, P_n^{alg}B_2^{alg}$ of $P^{alg}$ respectively
 by right $B_0^{alg}$, $B_1^{alg}$ and $B_2^{alg}$-modules. The map
\eqref{bruhat} is clearly filtration preserving so we get canonical maps of graded $k$-vector spaces
\[P_n^{alg}\otimes_{B_0^{alg}} k \to P_n^{alg}B_1^{alg}\otimes_{B_1^{alg}} k \oplus 
P_n^{alg}B_2^{alg}\otimes_{B_2^{alg}} k \]
It follows immediately from Theorem \ref{normalformalg} together with \eqref{secondfiltr}
that these maps induce isomorphisms on the associated graded vector space of the filtrations
(except in filtration 0 where the cokernel is $k$). This completes the proof.
\end{proof}

\begin{remark}
It would be interesting to know whether Theorem \ref{main2} holds under the assumption
that $H_*(B_i;k)$ are $H_*(B_0;k)$-flat but the maps $H_*(f_i;k)$ are not necessarily injective.
 This is true for discrete groups according to \cite[Theorem 4.1]{Fi}.
\end{remark}

\section{The homotopy decomposition of $G_\lambda$ for $1<\lambda\leq 2$.}

The following is a reformulation of \cite[Theorem 1.2, Theorem 3.1]{An} in a form which is suitable for
our purposes.
\begin{thm}[Anjos]
\label{silvia}
Let $k$ be a field. The diagram 
\[\xymatrix{
H_*(SO(3);k) \ar[r]^{\Delta} \ar[d] & H_*(SO(3)\times SO(3);k) \ar[d] \\
H_*(S^1\times SO(3);k) \ar[r] & H_*(G_\lambda;k)}\]
is a pushout square in the category of $k$-algebras.
\end{thm}
We are now ready to prove our main result.
\begin{proof}[Proof of Theorem \ref{main}]
The hypotheses of Theorem \ref{main1} are clearly satisfied. Therefore Theorem \ref{silvia}  
implies that the canonical map 
\[ (SO(3)\times SO(3)) \amalg_{SO(3)} (S^1 \times SO(3)) \lra G_\lambda \]
is a homology equivalence with any field coefficients. Hence it is an
equivalence on integral homology. Since both spaces are $H$-spaces it follows
that it is in fact a weak equivalence. This proves (i).

It remains to show the equivalence of (i)-(iii). 

(iii) $\Rightarrow$ (ii): This follows immediately by applying the Borel construction
\[EG_{\lambda} \times_{G_\lambda} - \] to the the weak equivalence of (iii).

(ii)  $\Rightarrow$ (i): By Theorem \ref{main2}, (ii) implies that the map $BP \to BG_\lambda$ is a weak
equivalence and looping it get that 
\[  P \lra G_\lambda \] 
is a weak equivalence.

(i)  $\Rightarrow$ (iii): Saying that the map (iii) is a weak equivalence amounts to
saying that 
\[ \hocolim (G_\lambda/(SO(3)\times SO(3)) \lla G_\lambda/SO(3) \lra G_{\lambda}/(SO(3)\times S^1) ) \]
is weakly contractible. Assuming (i), this follows from \eqref{eq1} and \eqref{eq2}.
\end{proof}

Here is an interesting consequence of Theorem \ref{main}(ii).
\begin{prop}
\label{fibersequence}
There is a fiber sequence
\begin{equation}
\label{fiberseq}
\Sigma^2SO(3) \lra BG_\lambda \lra BS^1\times BSO(3) \times BSO(3) 
\end{equation}
\end{prop}
\begin{proof}
Let $X=BSO(3)$ and $Y=BS^1$. Apply Theorem \ref{puppe} to the bottom half of the diagram,
where all vertical sequences are homotopy fiber sequences
\[ 
\xymatrix{ 
\Omega Y \ar[d] & \Omega X \times \Omega Y \ar[l] \ar[d] \ar[r] & \Omega X \ar[d] \\
X\times X \ar[d] & X \ar[l] \ar[d] \ar[r] & X \times Y \ar[d] \\
X\times X \times Y & X\times X \times Y \ar[l] \ar[r] & X\times X\times Y 
}\]
we see that the fiber of the map 
\[ BG_\lambda \to BS^1\times BSO(3) \times BSO(3) \]
is 
\[ \hocolim ( \Omega Y \lla \Omega X \times \Omega Y \lra \Omega X ) = \Omega X * \Omega Y 
\simeq \Sigma \left(\Omega X \wedge \Omega Y \right)\]
(where $\ast$ denotes the join). This completes the proof.
\end{proof}

In particular, we see that $G_\lambda$ has a canonical homotopy representation
\begin{equation}
\label{homap}
BG_\lambda \to BS^1 \times BSO(3) \times BSO(3) 
\end{equation}
even though, by Banyaga's theorem \cite[Theorem 10.25]{MS}, $G_\lambda$ is a simple group
and hence admits no nontrivial homomorphisms to compact Lie groups.

\begin{remark}
The identification of the homotopy type of $G_\lambda$ \cite[Theorem 1.1]{An}
is an immediate consequence of Proposition \ref{fibersequence}. In fact, looping \eqref{fiberseq}
we get a fiber sequence of loop maps
\[ 
\Omega \Sigma^2 SO(3) \lra G_\lambda \lra S^1 \times SO(3) \times SO(3) 
\] 
This has a section (given by including the subgroups in $G_\lambda$ and multiplying)
so $G_\lambda$ is weakly equivalent to the product of the base and the fiber.
\end{remark}

\begin{remark}
Let $X \subset SO(3)\times SO(3)$ denote the image of any section of the principal fiber bundle
\[ SO(3)\times SO(3) \lra SO(3) \]
determined by quotienting $SO(3)\times SO(3)$ by the diagonal $SO(3)$ action.
Denoting by $P$ the pushout of \eqref{G-diag}, it is not hard to see 
that the kernel of the canonical group homomorphism
\[ P \to S^1\times SO(3) \times SO(3) \]
is the free group generated by the set of commutators $[S^1,X] \subset P$ (this set is
independent of the choice of section). The space $[S^1,X] \subset P$  is homeomorphic 
to $\Sigma SO(3)$. Writing $FY$ for the free topological group
generated by the space $Y$, we have $FY \simeq \Omega \Sigma Y$. Thus \eqref{fiberseq} is 
the fiber sequence obtained by applying the classifying space functor to the short exact sequence
 of topological groups
\[ F[S^1,X] \lra P \lra S^1\times SO(3) \times SO(3). \]
\end{remark}

Using Theorem \ref{main}(ii) , it is easy to compute the ring of integral characteristic classes of 
symplectic fibrations with fibre $M_\lambda$.
\begin{cor}
\label{homology}
\[ H^*(BG_\lambda;\Z) = \Z[T,X_1,X_2,Y_1,Y_2,Z]/<T(X_i-Y_i),TZ,2X_1,2Y_1,2Z,Z^2>\]
where $|Y_i|=|X_i|=i+2$, $|T|=2$, and $|Z|=5$.
\end{cor}
\begin{proof}
Since 
\[  H^*(BSO(3)\times BSO(3);\Z) = \Z[X_1,X_2,Y_1,Y_2,Z]/<2X_1,2Y_1,2Z,Z^2> \]
this follows easily from the Mayer-Vietoris sequence
\[ \cdots \lra H^*(BG_\lambda ;\Z) \lra H^*(BSO(3)\times BSO(3);\Z) \oplus H^*(BS^1\times BSO(3);\Z) \]
\[ \llra{\phi} H^*(BSO(3);\Z)  \lra \cdots \]
which identifies $H^*(BG_\lambda;\Z)$ with the kernel of the map $\phi$. 
\end{proof}
In particular, the map $H^*(BS^1 \times BSO(3) \times BSO(3);\Z) \to H^*(BG_\lambda;\Z)$ 
induced by \eqref{homap} is surjective.

\begin{remark}
Using the computation of $H^*(BG_\lambda;\Q)$ by Abreu and McDuff \cite{AM}, 
Januszkiewicz and Kedra \cite{JK} have shown that all the characteristic classes in 
$H^*(BG_\lambda;\R)$ come from integrating monomials on the Chern classes of the vertical 
tangent bundle and the coupling class\footnote{ The coupling class $\Omega \in H^2( M_{hG_\lambda};\R)$ 
is determined by the requirements that it restricts to $[\omega_\lambda]$ on the fiber $M_\lambda$ 
and satisfies $\pi_*^{G_\lambda}(\Omega^{n+1})=0$ (cf. \cite[Chapter 6]{MS}).} over the fibers. 

In more detail, the fibration 
\begin{equation}
\label{fibint}
 M \hookrightarrow M_{hG_\lambda}:=M\times_{G_\lambda} EG_\lambda 
\rightarrow BG_\lambda
\end{equation}
yields a vertical tangent bundle 
$$T M_{hG_\lambda} := TM \times_{G_\lambda} EG_\lambda \rightarrow M_{hG_\lambda}.$$
Since $\cat J_\lambda$ is contractible, this bundle has (up to homotopy) a canonical complex structure.
Integrating monomials on the Chern classes $c_k$ and the coupling class 
$\Omega$ over the fiber of \eqref{fibint} yields elements in $H^*(BG_\lambda;\R)$.

If  $\alpha:T \rightarrow G_\lambda$ is a torus action on $M_\lambda$, the commutativity of 
the diagram
\[ \xymatrix{H^*(M_{hT})   \ar[d]^{\pi_*^T} & H^*(M_{hG_\lambda}) \ar[d]^{\pi_*^{G_\lambda}}  
\ar[l]_{\tilde{B}_\alpha^*} \\
H^{*-2n}(BT) & H^{*-2n}(BG_\lambda) \ar[l]_{B_\alpha^*}
}\]
together with the fact that $\pi_*^T$  can be computed by localization methods
gives information on the corresponding classes in $H^*(BG_\lambda)$.

For example, Januszkiewicz and Kedra's calculations \cite[Proposition 4.1.1]{JK} say
that the obvious map $BG_\lambda \lra BS^1$  derived from \eqref{homap} classifies the class 
$-\pi_*^{G_\lambda}(c_1^3)/8 \in H^2(BG_\lambda;\Z)$.
\end{remark}

\subsection{Acknowledgments} The first author would like to thank Joe Coffey for useful conversations.
The second author would like to thank Jer\^ome Scherer for pointing out    
that Whitehead's theorem is an immediate consequence of Puppe's theorem (cf. proof of 
Theorem \ref{main2}), and also Kasper Andersen and Jesper Grodal for useful discussions.

\end{document}